\documentclass[multi]{cambridge7A}
\usepackage[english,UKenglish]{babel}
\usepackage[longnamesfirst,sectionbib]{natbib}
\usepackage{chapterbib}
\usepackage{amsmath,amssymb,amsthm,amsfonts}
\usepackage{enumerate,url}
\usepackage{pict2e}

\theoremstyle{plain}
\newtheorem{theorem}{Theorem}[section]

\allowdisplaybreaks[1]

\setcitestyle{numbers,square,comma}

 \newcommand{\Reals}{{\mathbb{R}}}
  \newcommand{\Nats}{{\mathbb{N}}}
  \newcommand{\GG}{\mbox{${\mathcal G}$}}
    
 \newcommand{\PP}{\mbox{${\mathcal P}$}}
  
  \newcommand{\NN}{\mbox{${\mathbb{N}}$}}
  \newcommand{\TT}{\mbox{${\mathcal{T}}$}}
   \renewcommand{\SS}{\mbox{${\mathcal S}$}}
 \newcommand{\UP}{\NN_{(2)}}
  \newcommand{\ed}{\ \stackrel{d}{=} \ }
 \newcommand{\cd}{\ \stackrel{d}{\rightarrow} \ }
\newcommand{\bX}{\mbox{${\mathbf{X}}$}}
\newcommand{\bx}{\mbox{${\mathbf{x}}$}}
\newcommand{\bF}{\mbox{${\mathbf{F}}$}}
\newcommand{\bp}{\mbox{${\mathbf{p}}$}}

 \newcommand{\sfrac}[2]{{\textstyle\frac{#1}{#2}}}
  \newcommand{\Ex}{{\mathbb{E}}}
 \renewcommand{\Pr}{{\mathbb{P}}}

\providecommand*\Index[1]{#1\index{#1}}
\providecommand*\undex[1]{} 

\begin{document}
\alphafootnotes
\author[David J. Aldous]{David J. Aldous\footnotemark }
\chapter[More uses of exchangeability]{More uses of exchangeability:
  representations of complex random structures}
\footnotetext{Department of Statistics,
  Evans Hall, University of California at Berkeley, Berkeley, CA 94720-3860;
  aldous@stat.berkeley.edu; \url{http://www.stat.berkeley.edu/users/aldous}.
  Research supported by N.S.F. Grant DMS-0704159.}
\arabicfootnotes
\contributor{David J. Aldous
  \affiliation{University of California at Berkeley}}
\renewcommand\thesection{\arabic{section}}
\numberwithin{equation}{section}
\renewcommand\theequation{\thesection.\arabic{equation}}
\numberwithin{figure}{section}
\renewcommand\thefigure{\thesection.\arabic{figure}}

 \begin{abstract}
We review old and new uses of exchangeability, emphasizing the general theme of 
exchangeable representations of complex random structures.
Illustrations of this theme include 
processes of stochastic coalescence and fragmentation;
continuum random trees;  
second-order limits of distances in random graphs;
isometry classes of metric spaces with probability measures; 
limits of dense random graphs; and more sophisticated uses in finitary combinatorics.
 \end{abstract}
 
\subparagraph{AMS subject classification (MSC2010)}60G09, 60C05, 05C80


 \section{Introduction}
 \label{sec-INT}
\index{exchangeability|(}Kingman's\index{Kingman, J. F. C.} write-up \cite{MR0494344} of his 1977 Wald Lectures drew attention to the
subject of exchangeability,  and further indication of the topics of interest around
that time can be seen in the write-up \cite{me22} of my 1983 Saint-Flour lectures. 
As with any mathematical subject, one might expect some topics subsequently to
wither, some to blossom and new topics to emerge unanticipated.   This Festschrift paper
aims, in  informal lecture style,
\begin{enumerate}[(a)]
\item to recall the state of affairs 25 years ago (sections
\ref{sec-dFT}--\ref{sec-PEA},  \ref{sec-partitions});
\item to briefly describe three directions of subsequent development  that have
recently featured prominently in monographs
\cite{MR2253162,MR2161313,MR2245368} (sections \ref{sec-ff},
\ref{sec-partitions}--\ref{sec-bertoin});
\item to describe 3 recent rediscoveries, motivated by new theoretical topics
outside mainstream mathematical \Index{probability}, of the theory of representations of
\index{exchangeability!partial exchangeability}partially exchangeable arrays (sections
\ref{sec-isometry}, \ref{sec-dense}--\ref{sec-austin});
\item to emphasize a  general program that has interested me for 20 years.  It doesn't
have a standard name, but let me here call it \emph{exchangeable representations of
complex random structures}\index{exchangeability!exchangeable representation} (section \ref{sec-BP}).
\end{enumerate}

The survey focusses on mathematical probability; although the word 
\emph{Bayesian} appears several times, I do not attempt to cover the vast territory 
of explicit or implicit uses of exchangeability in Bayesian
statistics\index{Bayes, T.!Bayesian statistics|(}, except to mention here 
its use in \emph{hierarchical models}\index{Bayes, T.!Bayesian hierarchical model} \cite{LDA2003,MR2027492}.
 
This article is very much a bird's-eye view.
Of the monographs mentioned above, let me point out  Pitman's\index{Pitman, J. [Pitman, J. W.]}
\emph{Combinatorial Stochastic Processes} \cite{MR2245368}, which packs an extraordinary number of detailed results into 
200 pages of text and exercises.  
Exchangeability is a recurring theme in \cite{MR2245368}, which covers about
half of the topics we shall mention (and much more, not related to exchangeability), and so \cite{MR2245368}
is a natural starting place for the reader wishing to get to grips with details.

\section{Exchangeability}
\subsection{de Finetti's theorem}
\label{sec-dFT}
I use \emph{exchangeability} to mean, roughly, `applications of extensions of de
Finetti's\index{Finetti, B. de!de Finetti's theorem|(} theorem'.
Let me assume only that the reader is familiar with the 
definition of an exchangeable sequence of random variables
\[ (Z_i, i \geq 1) \ed (Z_{\pi(i)}, i \geq 1) \mbox{ for each finite permutation }
\pi \] 
and with the 
common verbal statement of de Finetti's theorem:
\begin{quote}
An infinite exchangeable sequence is distributed as a mixture of i.i.d. sequences. 
\end{quote}
Scanning the graduate-level textbooks on measure-theoretic probability on my
bookshelves, the theorem makes no appearance in about half, 
a  brief appearance in others  \cite{MR1324786,dur91v3,MR0651018}
and only three \cite{MR1476912,MR1422917,kall97} devote much more than a page to the
topic.

A remarkable feature of  de Finetti's theorem  is that there are many ways to state
essentially the same result, depending on the desired emphasis. This feature is best
seen when you state results more in words rather than symbols, so that's what I
shall do. 
Take a nice space $S$, and either don't worry what `nice' means, or assume $S =
\Reals$. 
Write $\PP(S)$ for the space of probability measures on $S$.  
Write $\mu$ for a typical element of $\PP(S)$ and write $\alpha$ for a typical {\em
random} element of $\PP(S)$, that is a typical \emph{random
measure}\index{random measure}.  
When we define an infinite exchangeable sequence of $S$-valued random variables we are really defining an 
\emph{exchangeable measure} ($\Theta$, say) on $\PP(S^\infty)$, where $\Theta$ is the distribution of the sequence.

\paragraph{Functional analysis viewpoint.} 
The subset 
$\{\mu^\infty = \mu \times \mu \times \mu \times \ldots : \mu \in \PP(S)\} \subset
\PP(S^\infty)$ 
is the set of extreme points of the  convex set of exchangeable elements of
$\PP(S^\infty)$, and the identification
\[ \Theta(\cdot) = \int_{\PP(S)} \mu^\infty(\cdot) \ \Lambda(d \mu) \]  gives a
bijection between probability measures $\Lambda$ \underline{on} $\PP(S)$ 
(that is, $\Lambda \in \PP(\PP(S))$) and exchangeable measures $\Theta$ in
$\PP(S^\infty)$.

\paragraph{Probability viewpoint.} 
Here are successively more explicit versions of the same idea.
Let $(Z_i, \ 1 \leq i <  \infty)$ be exchangeable $S$-valued.
\begin{enumerate}[(a)]
\item Conditional on the tail (or invariant or exchangeable) $\sigma$-field of the
sequence $(Z_i)$, the random variables $Z_i$ are i.i.d.
\item There exists a random measure $\alpha$ such that, conditional on $\alpha = \mu$,
the random variables $Z_i$ are i.i.d.($\mu$).
\item  Giving $\PP(S)$ the topology of
\undex{weak convergence}weak convergence, the empirical measure $F_n =
F_n(\omega, \cdot) = n^{-1} \sum_{i=1}^n 1_{(Z_i(\omega) \in \cdot)}$ 
converges a.s. to a limit random measure $\alpha(\omega, \cdot)$
satisfying (b).
\end{enumerate}

\paragraph{Theoretical statistics viewpoint.} 
In contexts where a frequentist would assume data are i.i.d.($\mu$)  from an unknown
distribution $\mu$, a Bayes\-ian would put a prior distribution  $\Lambda$ on possible
$\mu$; so de Finetti's theorem is saying that the Bayesian assumption is logically
equivalent to the assumption that the data 
$(Z_i, i \geq 1)$ are exchangeable.  Note a mathematical consequence.  There is a
posterior distribution 
$\Lambda_n(\omega, \cdot) \in  \PP(\PP(S))$ for $\Lambda$ given
$(Z_1,\ldots,Z_n)$, and an extension of (c) above is
\begin{enumerate}[(a)]\addtocounter{enumi}{3}
\item $\Lambda_n(\omega, \cdot) \to \delta_{\alpha(\omega, \cdot)}$ a.s.  in
$\PP(\PP(S))$.
\end{enumerate}
Such results are historically often used as a starting point for
philosophical and mathematical discussion of consistency/inconsistency  of
frequentist and Bayesian methods,
inspired of course by Bruno de Finetti\index{Finetti, B. de} himself\index{Bayes, T.!Bayesian statistics|)}.

\vspace{0.2in}
\noindent
But there's more!  
Later we encounter at least two further, somewhat different, viewpoints:  explicit 
constructions (section \ref{sec-constructdF}),  and our central theme of using exchangeability to describe 
complex structures (sections \ref{sec-Using} and \ref{sec-BP}).  
This theme is related to the general features that
\begin{itemize}
\item exchangeable-like properties are preserved under weak convergence;
\item parallel to representation theorems for infinite exchangeable-like structures, are convergence theorems giving necessary and sufficient condition for finite exchangeable-like structures to converge in distribution to an infinite such structure.
\end{itemize}
In the setting of de Finetti's theorem, the condition for finite exchangeable sequences 
$\bX^{(n)} = (X^{(n)}_1,\ldots, X^{(n)}_n)$ to converge in distribution to an infinite exchangeable sequence $\bX$ is 
\[ \alpha_n \cd \alpha \mbox{ on }  \PP(S) \]
where $\alpha$ is the `directing' random measure for $\bX$  in (b) above, and $\alpha_n$ is the empirical distribution of  $(X^{(n)}_1,\ldots, X^{(n)}_n)$.
Note that when we talk of convergence in distribution to infinite sequences or arrays, we mean w.r.t.  product topology, i.e. convergence of finite restrictions.

\subsection{Exchangeability, 25 years ago}
\label{sec-25y}
Here I list topics from the two old surveys  \cite{MR0494344,me22}, for the purpose
of  saying a few words about those topics I will \emph{not} mention further, 
while pointing to sections where other topics will be discussed further.  

\paragraph{Classical topics not using de Finetti's theorem.}
\begin{enumerate}[(a)]
\item Combinatorial aspects  for classical stochastic processes, e.g. \index{ballot theorem}ballot theorems:
 \cite{MR0217858}.
\item Weak convergence for \index{sampling|(}`sampling without replacement' processes (e.g.
\cite{MR0233396}  Thm 24.1). 
\end{enumerate}

\paragraph{Variants of de Finetti's theorem.} 
Several variants were already classical, for instance:
\begin{enumerate}[(a)]\addtocounter{enumi}{2}
\item \index{Schoenberg, I. J.}\index{Diaconis, P.}\index{Freedman, D. A.}Schoenberg's\footnote{Persi Diaconis observes that the result is hard to deduce from Schoenberg 
\cite{MR1501980} and should really be attributed to Freedman \cite{MR0156369}.}
 theorem (\cite{me22}(3.6))   for the special case of spherically
symmetric sequences;
\item the analogous representation (\cite{me22}(3.9)) in the setting of two
se-\break quences 
$(X_i, 1 \leq i < \infty; \ Y_j, 1 \leq j < \infty)$ whose joint
distribution is invariant under finite permutations of either;
\item the \emph{selection property} (\cite{MR0494344} p.188), that the exchangeability
hypothesis in  de Finetti's theorem can be weakened to the assumption
\begin{align*}
 (X_1,X_2,\ldots X_n) &\ed (X_{k_1}, X_{k_2},\ldots,X_{k_n})\\
 &\mbox{ for all } 
 1 \leq k_1 < k_2 < \cdots < k_n .
\end{align*}
\end{enumerate}

Other variants had been developed in the 1970s, for instance:
\begin{enumerate}[(a)]\addtocounter{enumi}{5}
\item the analog for continuous-time processes with
exchangeable\index{exchangeability!exchangeable increments} \emph{increments}
\cite{MR0394842};
\item Kingman's\index{Kingman, J. F. C.!influence}
\index{Kingman, J. F. C.!Kingman paintbox}paintbox theorem for \index{random partition}exchangeable random partitions; see section
\ref{sec-partitions}.
\end{enumerate}

\paragraph{Finite versions.}
The general forms of  de Finetti's theorem and some classical variants can be proved
by comparing \index{sampling|)}sampling with and without replacement. This method \cite{MR577313} also 
yields finite-$n$
variants\index{Finetti, B. de!de Finetti's theorem|)}.

\paragraph{Mathematical population genetics, the \index{coalescence}coalescent and the
Pois\-son--Dirichlet distribution.}\index{mathematical
genetics}\index{Poisson, S. D.!PoissonDirichlet distribution@Poisson--Dirichlet distribution}
Exchangeability is involved in this  large circle of ideas, developed in part by
Kingman in the 1970s, which continues to prove fruitful in many ways.
For the population genetics aspects of Kingman's work see the article by
Ewens\index{Ewens, W. J.} and
Watterson\index{Watterson, G. A.} \cite{EW10} in this volume;  also the \emph{Kingman coalescent} fits into the more general \emph{stochastic coalescent} material in section \ref{sec-bertoin}.

\paragraph{The \Index{subsequence principle}.}  The idea emerged in the 1970s that, from any
tight sequence of random variables, one can extract a subsequence which is close to
exchangeable,   close enough that one can prove analogs of classical limit theorems
(CLT and LIL, for instance) for the subsequence.  General versions of this principle
were established in  \cite{me2,MR859840}, which pretty much killed the topic.

\paragraph{Sufficient statistics and mixtures of Markov chains.} 
One can often make a Bayesian interpretation\index{Bayes, T.!Bayesian statistics|(} of `\Index{sufficient statistic}'  in terms of
some context-dependent invariance 
property  \cite{MR786142}.   Somewhat similarly, one can characterize mixtures of
\index{Markov, A. A.!Markov chain}Markov chains via the property that transition counts are sufficient statistics
\cite{MR556418}.

\subsection{Partially exchangeable arrays}
\label{sec-PEA}
The topic, emerging around 1980, of \emph{partially exchangeable arrays}\index{exchangeability!partial exchangeability|(}, plays a
role in what follows and so requires more attention. 
Take a measurable function 
$f: [0,1]^2 \to S$ which is symmetric, in the sense  
$f(x,y) \equiv  f(y,x)$. Take $(U_i, i \geq 1)$ i.i.d. Uniform$(0,1)$ and consider
the array 
\begin{equation}
 X_{\{i,j\}} := f(U_i,U_j)  \label{basic_array}
 \end{equation}
indexed by the set $\UP$ of unordered pairs $\{i,j\}$.  
The exchangeability property of $(U_i)$ implies what we shall call the  {\em
partially exchangeable} property for the array:
\begin{align}
( X_{\{i,j\}} , \ {\{i,j\}} \in \UP) &\ed  ( X_{\{\pi(i),\pi(j)\}} , \ {\{i,j\}} \in
\UP)\notag \\
  &\qquad\mbox{ for each finite permutation } \pi  .
\label{def-PE}
\end{align}
  Note this is a weaker property
than the `fully exchangeable' property for the countable collection 
$( X_{\{i,j\}} , \ {\{i,j\}} \in \UP)$, because the permutations of $\UP$ which are
of the particular form 
${\{i,j\}} \to {\{\pi(i),\pi(j)\}}$ for a finite permutation $\pi$ of $\NN$ are only
a subset of all permutations of $\UP$.  

Aside from construction (\ref{basic_array}), how else can one produce an array with
this  partially exchangeable property? 
Well, an array with i.i.d. entries has the property, and so does the trivial case
where all entries are the same r.v.  
After a moment's thought we realize we can combine these ideas as follows.

Take a function $f: [0,1]^4 \to S$ such that $f(u,u_1,u_2,u_{12})$ is symmetric in
$(u_1,u_2)$, and then define 
\begin{equation}
 X_{\{i,j\}} := f(U,U_i,U_j, U_{\{i,j\}}  )  
 \label{PE-rep}
 \end{equation}
 where all the r.v.s in the families $U, (U_i, i \in \NN), (U_{\{i,j\}} ,
\ {\{i,j\}} \in \UP)$ are i.i.d. Uniform$(0,1)$.
Then the resulting array $\bX = (X_{\{i,j\}} )$ is partially exchangeable.

Finding oneself unable to devise any other constructions, it becomes natural to conjecture
that every partially exchangeable array has a representation (in distribution) of
form (\ref{PE-rep}).
This was proved by Hoover\index{Hoover, D. N.} \cite{hoover-rel} and (in the parallel non-symmetric setting) by Aldous
\cite{me11}, the latter proof having been substantially
simplified due to a personal communication from Kingman\index{Kingman, J. F. C.}.

Constructions of partially exchangeable arrays appear in Bayesian
statistical\index{Bayes, T.!Bayesian statistics|)} modeling; 
see e.g. the family of copulae\index{copula} introduced in \cite{MR2026570} in the context of a 
semi-parametric model for \Index{Value at Risk}.

\subsection{Fast forward}
\label{sec-ff}
Such \emph{partially exchangeable representation theorems} were the state of the art
in the 1984 survey \cite{me22}.
They were subsequently  
extended systematically by Kallenberg\index{Kallenberg, O.}, both for arrays and analogs such as
exchangeable-increments continuous-parameter processes, and for the \emph{rotatable 
matri\-ces}\index{rotatable matrix} to be mentioned in section \ref{sec-matrices}, 
 during the late 1980s and early 1990s.   The whole topic of
representation theorems 
 is the subject of Chapters 7--9 of Kallenberg's 2005 monograph 
\cite{MR2161313}.  Not only does this monograph  provide a  canonical reference to
the theorems, 
but also its introduction  provides an excellent summary
 of the topic.  
 
 In the particular setting above we have
\begin{theorem}[Partially Exchangeable Representation Theorem]
\label{T1}
An array $\bX$ which is partially exchangeable, in the sense\/ $(\ref{def-PE})$, has a representation in the form\/ $(\ref{PE-rep})$.
\end{theorem}
This is one of the family of results described carefully in Chapter 7 of \cite{MR2161313}.  
There are analogous results for higher-parameter arrays $(X_{ijk})$, and for arrays in which the 
`joint exchangeability' assumption (\ref{def-PE}) is replaced by a `separate exchangeability' assumption for non-symmetric arrays
$(X_{i,j}, 1 \le i,j < \infty)$:
\begin{align*}
 (X_{i,j}, 1 \le i,j < \infty) &\ed (X_{\pi_1(i), \pi_2(j)}, 1 \le i,j < \infty)\\
 &\qquad\mbox{ for finite permutations } \pi_1, \pi_2  .
\end{align*}

One aspect of this theory is surprisingly subtle, and that is the 
issue of \emph{uniqueness} of representing functions $f$. 
In representation (\ref{PE-rep}), if we take Lebesgue-measure-preserving
maps\index{measure preserving map@measure-preserving map}
$\phi_0$, $\phi_1$, $\phi_2$ from $[0,1]$ to $[0,1]$, then the arrays $\bX$ and $\bX^*$ obtained from $f$ and from 
$f^* (u,u_1,u_2,u_{12}) :=  f(\phi_0(u), \phi_1(u_1),\phi_1(u_2),\phi_2(u_{12}))$ must have the same distribution.
But this is not the only way to make arrays have the same distribution: 
there are other ways to construct measure-preserving 
transformations of $[0,1]^4$, and 
(because measure-preserving transformations are not invertible in general)
one needs to insert randomization variables. 
(I thank the referee for correcting a blunder in my first draft, and for the comment 
``this may be a major reason why
\index{nonstandard analysis@non-standard analysis}non-standard analysis is effective here''.)
For an explicit statement of the uniqueness result in two dimensions see 
\cite{MR1003713} and for higher dimensions see \cite{MR2161313}.

Relative to the Big Picture of Mathematics, 
this theory of partial exchangeability\index{exchangeability!partial exchangeability|)} was perhaps regarded during 1980--2005 as a rather small niche inside mathematical \Index{probability}---and ignored outside mathematical probability.   So
it is ironic that around 2004--8 it was rediscovered in at least three different
contexts outside mainstream mathematical probability.  Let me say one such context right now and
the others later (sections \ref{sec-dense} and \ref{sec-austin}).

\subsection{Isometry classes of metric spaces with probability measures}  
\label{sec-isometry}\index{isometry}
The definition of \emph{isometry} between two metric spaces 
$(S_1,d_1)$ and $(S_2,\allowbreak d_2)$ contains an `if there exists \ldots'
expression.  Asking for a \emph{characterization} of metric spaces up to isometry is
asking for a scheme that associates some notion of `label' to each metric space in
such a way that two metric spaces are isometric if and only if they have the same
label.  
I am not an expert on this topic, but I believe there is no known such
characterization.  

But suppose instead we consider `metric spaces with probability measure',
$(S_1,d_1,\mu_1)$ and $(S_2,d_2,\mu_2)$, and require the isometry to map $\mu_1$ to
$\mu_2$.  It turns out there is now a conceptually simple
characterization.  Given
$(S,d,\mu)$, take i.i.d.($\mu)$ random elements  $(\xi_i, 1 \le i < \infty)$ of $S$,
form the array (of form (\ref{basic_array}))
\begin{equation}
 X_{\{i,j\}} = d(\xi_i, \xi_j);\  \{i,j\} \in \UP \label{Xdxi}
\end{equation}
and let $\Psi$ be the distribution of the infinite random array.   It is obvious
that, for two isometric 
`metric spaces with probability measure', we get the same $\Psi$, and the converse
is a simple albeit technical consequence of the uniqueness 
part of Theorem \ref{T1},
implying:
\begin{gather}
\mbox{`metric spaces with probability measure' are characterized}\notag
 \\[-1mm]
\mbox{up to isometry by the distribution $\Psi$.}\label{vershik-cor}
\end{gather}
This result was given by 
Vershik\index{Vershik, A. M.} \cite{MR2086637}, as one rediscovery of part of the general theory of
partial exchangeability.

\subsection{Rotatable arrays and random matrices with symmetry properties}
\label{sec-matrices}\index{random matrix}\index{rotatable array}
In Theorem \ref{T1} we  described $\bX = (X_{ij})$ as an \emph{array} instead of a \emph{matrix}, partly because of the extension to higher-dimensional parametrizations and partly because we never engage matrix multiplication.  
Now regarding $\bX$ as a matrix, one can impose stronger `matrix-theoretic' assumptions  
and ask for characterizations of the random matrices satisfying such assumptions.
One basic case, \emph{rotatable matrices}, is where the $n \times n$ restrictions are invariant in distribution under the orthogonal group, and the characterization goes back to 
\cite{me11}.  Two other cases (I thank the referee for suggesting (ii) and the
subsequent remark) are
\begin{enumerate}[(i)]
\item non-negative definite jointly exchangeable arrays: \cite{MR666087,panch09};
\item rotatable \emph{Hermitian} matrices\index{Hermite, C.!Hermitian matrix} \cite{MR1402920}, motivated indirectly 
by problems in
\index{quantum!mechanics}quantum mechanics and thereby related to the huge
literature on semicircular laws\index{semicircle distribution}.  
\end{enumerate}

Returning to the basic case of rotatable matrices, for the higher-dimensional analogs 
 the basic representations are
naturally stated in terms of \index{multiple Wiener--It{\^o} integral}multiple Wiener--It{\^o} integrals, which form
the fundamental examples of rotatable random functionals\index{rotatable random functional}.  Such multiple Wiener--It{\^o}
integrals are also a basic tool in \index{Malliavin, P.!Malliavin calculus} \cite{MR2200233}, 
a subject with important applications to analysis.

\subsection{Revisiting de Finetti's theorem}
\label{sec-constructdF}\index{Finetti, B. de!de Finetti's theorem}
Returning to a previous comment, the theory of  partially exchangeable
representation theorems reminds us that one can take a similar view of de Finetti's
theorem itself, 
to add to the list in section \ref{sec-dFT}.

\paragraph{Construction viewpoint.} 
Given a measurable function $f: [0,1]^2 \to S$ and i.i.d. Uniform$(0,1)$ random
variables 
$(U; U_i, i \geq 1)$, the process 
$(Z_i, i \geq 1)$ defined by 
$Z_i = f(U,U_i)$ is exchangeable, and every exchangeable process arises (in
distribution) in this way from some $f$.

\section{Using exchangeability to describe complex structures}
\label{sec-Using}
Here is my attempt at articulating the first part of the central theme of this paper\index{exchangeability!exchangeable representation}.
\begin{quote}
One way of examining a complex mathematical structure is to sam\-ple i.i.d. random
points and look at some 
form of induced substructure relating the random points. 
\end{quote} 
The idea being that the i.i.d. sampling induces some kind of
`exchangeability' on the distribution of the substructure, when the substructure is
regarded as an object in its own right.

The `\Index{isometry}' result (\ref{vershik-cor})  nicely fits this theme---the
substructure is simply the induced metric on the sampled points. The rest of the
present paper seeks to illustrate that this, admittedly very vague, way of looking
at structures
can indeed be useful, conceptually and/or technically.  Let us mention here two
prototypical examples (which will reappear later) of what a `substructure' might
be.  
Given $k$ vertices $v(1)$, \ldots, $v(k)$ in a \Index{graph}, one
can immediately see two different ways to define an induced substructure.
\begin{enumerate}[(i)]
\item The induced subgraph on vertices 1, \ldots, $k$: there is an edge $(i,j)$ iff the
original graph has an edge $(v(i),v(j))$.
\item The distance matrix: $d(i,j)$ is the number of edges in the shortest path from
$v(i)$ to $v(j)$.
\end{enumerate}
But before considering graph theoretic examples, let us explain with hindsight how
Kingman's\undex{Kingman, J. F. C.|(} work 
on \index{random partition|(}exchangeable random partitions fits this theme.

\subsection{Exchangeable random partitions and Kingman's paintbox theorem}\index{Kingman, J. F. C.!Kingman paintbox|(}
\label{sec-partitions}
The material here is covered in detail in Pitman\index{Pitman, J. [Pitman, J. W.]} \cite{MR2245368} Chapters 2--4.

Given a discrete sub-probability distribution, one can write the probabilities in
decreasing order as 
$p_1 \geq p_2 \ge \ldots > 0$ 
and then write $p_{(\infty)}:= 1 - \sum_j p_j \ge 0$ to define a probability
distribution $\bp$.
Imagine objects 1, 2, 3, \ldots\ each independently being colored, assigned color
$j$ with probability 
$p_j$ or assigned with probability $p_{(\infty)}$ a unique color  
(different from that assigned to any other object).  Then consider the
resulting `same color' equivalence 
classes as a random partition of $\Nats$.    So a realization of this process might be
\begin{equation}
 \{1,5,6,9,13, \ldots\}, \ \ \{2,3,8,11, 15,\ldots\}, \ \ \{4 \}, \ \ \{7,23, \ldots \}, \
\ldots \ldots  ;
\label{sbo}
\end{equation}
sets in the partition are either infinite or singletons.
This \emph{paintbox}($\bp$) distribution on partitions is exchangeable in the natural
sense.  
\emph{Kingman's paintbox theorem}, an analog of de Finetti's theorem, states that
every exchangeable 
random partition of $\Nats$ 
is distributed as a mixture over $\bp$ of paintbox($\bp$) distributions.

We mentioned in section \ref{sec-dFT} as a general feature that, associated with a representation theorem like this, 
there will be a convergence theorem.  Here are two slightly different ways of looking at the convergence theorem in the present setting.
Suppose that for each $n$ we are given an arbitrary random (or non-random) 
partition $\GG(n)$ of $\{1,2,\ldots,n\}$.  
For each $k < n$ sample without replacement $k$ times from $\{1,2,\ldots,n\}$ to get 
$U_n(1)$, \ldots, $U_n(k)$, consider the induced partition on the sampled elements
$U_n(1)$, \ldots, $U_n(k)$, 
and relabel these elements as 1, \ldots, $k$ to get a random partition $\SS(n,k)$ of
$\{1,2,\ldots,k\}$. 
This random partition $\SS(n,k)$ is clearly exchangeable.  If there is a limit 
\begin{equation}
 \SS(n,k) \cd \SS_k \mbox{ as } n \to \infty 
\label{part-ss}
\end{equation}
(the set of all possible partitions of $\{1,2,\ldots,k\}$ is finite, 
so there is nothing technically sophisticated here) 
then the limit $\SS_k$ is exchangeable; and if a limit (\ref{part-ss}) exists for
all $k$ 
then the family $(\SS_k, 1 \le k < \infty)$ specifies the distribution of an 
exchangeable random partition of $\Nats$, to which Kingman's paintbox theorem can be
applied.

The specific phrasing above was chosen to fit a general framework in section \ref{sec-program} later, 
but here is a more natural phrasing.  
For any random partition of $\{1,\ldots,n\}$ write 
$\bF^{(n)} = (F^{(n)}_1, F^{(n)}_2, \ldots)$ for the \emph{ranked empirical frequencies}\index{ranked empirical frequencies}, the numbers 
$n^{-1} \times $(sizes of sets in partition) in decreasing order.  
For a paintbox($\bp$) distribution the SLLN\index{strong law of large numbers (SLLN)} implies 
$\bF^{(n)} \to \bp$ a.s., and so Kingman's paintbox theorem implies that for any infinite exchangeable random partition $\Pi$, 
the limit $\bF^{(n)} \to \bF$ exists a.s. and is the  `directing random measure' 
(conditional on $\bF = \bp$ the distribution of $\Pi$ is paintbox($\bp$)).
Now suppose for each $n$ we have an \emph{exchangeable} random partition $\Pi^{(n)}$ of $\{1,2,\ldots,n\}$ and write $\bF^{(n)}$ for its ranked empirical frequencies.
The \emph{convergence theorem} states that the sequence $\Pi^{(n)}$ converges in distribution 
(meaning its restriction to $\{1,\ldots,k\}$ converges, for each $k$) to some limit $\Pi$, which is necessarily some infinite exchangeable random partition 
with some directing random measure $\bF$, if and only if 
$\bF^{(n)} = (F^{(n)}_j, \  1 \leq j < \infty)  \cd \bF  = (F_j, \  1 \leq j < \infty)  $. 

A final important idea is \emph{size-biased order}\index{size biasing}.
In the context of exchangeable random partitions this just means writing the 
components in a realization, as at (\ref{sbo}), starting with the component
containing element $1$,
then the component containing the least element not in the first component, and so
on.  
In the infinite case, the frequencies 
$\bF^* = (F^*_1, F^*_2, \ldots)$ of the components in size-biased order are just a
random permutation of 
the frequencies $\bF$ given by Kingman's paintbox theorem.  
In the paintbox($\bp$) case, replacing non-random $\bF = \bp$ by random $\bF^*$ is
perhaps merely
complicating matters, but in the general case of random $\bF$ it is often more
natural to 
work with the size-biased order than with the ranked order. 
For instance, the size-biased order codes information such as
\[ \Ex (F^*_1)^m = \Ex \sum_{i \geq 1} F_i^{m+1} . \]

I am highlighting these `structural' results as part of my overall theme, 
but in many ways the concrete examples are more interesting.  
The one-parameter  \index{Poisson, S. D.!PoissonDirichlet distribution@Poisson--Dirichlet distribution}Poisson--Dirichlet$(\theta)$ family was already recognized 25 years ago as a mathematically
canonical family of measures arising in several different contexts:
the \index{Ewens, W. J.!Ewens sampling formula (ESF)}Ewens sampling formula in neutral population
genetics\index{mathematical genetics},  
the `\index{Chinese restaurant process (CRP)}Chinese restaurant process' construction, 
a construction via \index{subordinator}subordinators, 
the size-biased order of asymptotic frequencies is the GEM distribution;
and special cases arise as limits of  component sizes in \index{random
permutation}random permutations and in  random mappings\index{random mapping}.
Subsequently the two-parameter Poisson--Dirichlet$(\alpha,\theta)$
distribution introduced by Pitman--Yor\index{Pitman, J. [Pitman, J. W.]}\index{Yor, M.} \cite{MR1434129} was shown to
possess many analogous properties.   
The paper \cite{GHP09} in this volume gives the flavor of current work in this direction.

Now partitions are rather simple structures, and the paintbox theorem\undex{Kingman, J. F. C.|)}\index{Kingman, J. F. C.!Kingman paintbox|)} (which
can be derived from  de Finetti's\index{Finetti, B. de!de
Finetti's theorem} theorem) isn't so convincing as an illustration of the theme `using exchangeability to describe complex structures'.  
The theme becomes more visible when we consider the more complex setting of partitions evolving over time,  and this setting arises naturally in the following context\index{random partition|)}.

\subsection{Stochastic coalescence and fragmentation}
\label{sec-bertoin}\index{coalescence|(}\index{fragmentation|(}
The topic of this section is treated in detail in Bertoin\index{Bertoin, J.}
\cite{MR2253162}, the third monograph in which exchangeability has
recently featured prominently.
The topic concerns models in which, at each time, unit mass is split into
clusters
of masses $\{x_j\}$.
One studies models of dynamics under which clusters split
(\emph{stochastic fragmentation}) or merge (\emph{stochastic coalescence} or
\emph{coagulation}\footnote{The word \emph{coagulation}, introduced in
German in \cite{smoluch}, sounds strange
to the native English speaker  to
 whom it suggests blood clotting; \emph{coalescence} seems a more apposite
English word.})
according to some random rules.

Conceptually, states are unordered collections $\{x_j\}$ of positive
numbers with
$\sum_j x_j = 1$.
What is a good mathematical representation of such states?
The first representation one might devise is as vectors
in decreasing order, say
$\bx^{\downarrow} = (x_1,x_2,\ldots )$.
But this representation has two related unsatisfactory features;
fragmenting one
cluster
forces one to relabel the others; and
given the realizations at two times, one can't
identify a particular cluster at the later time as a fragment of a particular
cluster
at the earlier time.
These difficulties go away if we think instead in terms of sampling
`atoms' and tracking which cluster they are in.
A uniform random atom will be in a mass-$x$ cluster of a configuration
$\bx^{\downarrow} $  with probability $x$;
sampling atoms $i = 1$, 2, 3, \ldots\ and taking the partition of
$\{1,2,3,\ldots\}$  into  `atoms of the same cluster' components
gives an exchangeable random partition $\Pi$ with
paintbox($\bx^{\downarrow} $) distribution.

Thus instead of representing a process as
$(\bX^{\downarrow}(t), 0 \leq t < \infty)$
we can represent it as a partition-valued process
$(\Pi(t),   0 \leq t < \infty)$ which tracks the positions of (i.e. the
clusters containing) particular  atoms.
For fixed $t$, both $\Pi(t)$ and $\bX^{\downarrow}(t)$ give the same
information about the cluster masses---and note
that clusters in $\Pi(t)$ automatically appear in size-biased order.
But as processes in $t$,
$(\Pi(t),   0 \leq t < \infty)$ gives more information than
$(\bX^{\downarrow}(t), 0 \leq t < \infty)$,
and in particular avoids the unsatisfactory features mentioned above.

Now in one sense this is merely a technical device, but I find it does give some
helpful insights.

\paragraph{The basic general stochastic models.}
In the basic model of  \emph{stochas\-tic fragmentation},
different  clusters evolve independently, a mass-$x$ clus\-ter splitting at
some
stochastic rate
$\lambda_x$ into clusters whose relative masses $(x_j/x, \ j \geq 1)$
follow some
probability distribution $\mu_x(\cdot)$.
(So the model neglects detailed $3$-dimensional geometry; the shape of a
cluster
is assumed not to affect its propensity to split, and different clusters
do not
interact).
Especially tractable is the \emph{self-similar}\index{self-similar} case where $\mu_x = \mu_1$ and
$\lambda_x = x^\alpha$ for some \emph{scaling exponent}\index{scaling exponent} $\alpha$.   Such
processes
are closely related to classical topics in theoretical and applied
probability---the log-masses form a continuous time branching random
walk\index{branching random walk (BRW)}, and the mass
of the
cluster containing a sample atom  forms a
continuous-time \index{Markov, A. A.!Markov process}Markov process on state space $(0,1]$.

The basic model for  \emph{stochastic coalescence} is to have $n$
particles, initially in
single particle clusters of masses $1/n$, and let clusters merge according
to a kernel
$\kappa(x,x^\prime)$ indicating the rate
(probability per unit time) at which a typical pair of clusters of masses $x$
and $x^\prime$ may merge.
For fixed $n$ this is just a finite-state continuous-time \index{Markov, A. A.!Markov chain}Markov chain,
but it is natural to study $n \to \infty$  limits, and there are two
different regimes.
On the time-scale where typical clusters contain $O(1)$ particles, i.e.
have mass $O(1/n)$,
there is an intuitively natural
\emph{hydrodyamical limit}\index{hydrodyamical limit} (law of large numbers), that is differential
equations for the relative proportions
$y_i(t)$ of $i$-particle clusters in the $n \to \infty$ limit.  This  {\em
Smoluchowski coagulation equation}\index{Smoluchowski, M. R. von
S.!Smoluchowski coagulation equation}
has a long history in several areas of science such as \Index{physical chemistry}, as
indicated in the survey \cite{me78}.
Recent theoretical work has made rigorous the connection between the
stochastic and
deterministic models, and part of this is described in \cite{MR2253162} 
Chapter 5.
A different limit regime concerns the time-scale when the largest clusters 
contain order $n$ particles, i.e. have mass of order $1$.  
In this limit we have real-valued cluster sizes evolving over time $(-\infty, \infty)$ 
and `\index{started from dust}starting with dust' at time $- \infty$, that is with the largest cluster
mass $\to 0$ as $t \to - \infty$ 
(just as, in the basic fragmentation model, the largest cluster
mass $\to 0$ as $t \to + \infty$ ) and $\to 1$ as $t \to + \infty$. 
(So these models incidentally provide novel examples within the classical 
topic of \emph{entrance boundaries}\index{entrance boundary} for Markov
processes\index{Markov, A. A.!Markov process|(}).

Finally recall \emph{Kingman's coalescent}\index{Kingman,
J. F. C.!influence}\index{Kingman, J. F. C.!Kingman coalescent|(}, as a model of
genealogical lines of descent\index{line of descent} within
neutral population genetics, which (with its many subsequent
variations) has become a recognized topic within mathematical population
genetics\index{mathematical genetics}---see e.g.  Wakeley  \cite{wakeley}. 
Although the background story is different, it can mathematically
be identified with the constant  rate ($\kappa(x,x^\prime) = 1$) stochastic
coalescent in the present context.

\paragraph{Discussion and special cases.}There are three settings above 
(fragmentation; discrete-particle coalescence; continuous-mass coalescence)
which one might formalize differently, but the advantage of the
`\index{random partition}exchangeable
random partition' set-up 
is that each can be represented as a partition-valued process $(\Pi(t))$.
Intuitively, coalescence and fragmentation are time-reversals\index{time reversal} of each
other, and it is noteworthy that
\begin{enumerate}[(i)]
\item there are several fascinating examples of special models where a
precise \Index{duality} relation exists and is useful (see e.g.\
section \ref{sec-compl} (iv));
\item but there seems to be no \emph{general} precise duality relationship within
the usual stochastic models.
\end{enumerate}

In the general models the processes $(\Pi(t))$ are all Markov, as processes on
partitions of
$\Nats$.
One can consider their restrictions $(\Pi_k(t))$ to partitions of
$\{1,\ldots,k\}$,
i.e. consider masses of the components containing $k$ sampled atoms.
In general $(\Pi_k(t))$ will not be Markov, but it is
Markov in special cases, which are therefore particularly
tractable. 
One such case (\cite{MR2253162} section 3.1) is \emph{homogeneous
fragmentation}, where each cluster has the same splitting rate
($\lambda_x \equiv \lambda_1$).
Another such case (\cite{MR2253162} Chapter 4)
is the elegant general theory of
\emph{exchangeable coalescents}, which
eliminates the
`only binary merging' aspect of Kingman's coalescent, and is
interpretable as $n \to \infty$ limit genealogies
of more general models in population
genetics\index{coalescence|)}\index{fragmentation|)}\index{Kingman,
J. F. C.!Kingman coalescent|)}.

\section{Construction of, and convergence to, infinite random
combinatorial objects}
\label{sec-BP}
\subsection{A general program}
\label{sec-program}
de Finetti's\index{Finetti, B. de!de Finetti's theorem} theorem refers specifically to \emph{infinite} sequences. Of course we
can always try to view an infinite object as a limit of finite objects, and in the
`25 years ago' surveys such convergence ideas were explicit in
the context of weak convergence for `\Index{sampling} without replacement' processes, 
in finite versions such as \cite{MR577313}, 
and in some other contexts, such as Kingman's theory of exchangeable random
partitions.
I previously stated the first part of our central theme as 
\begin{quote}
One way of examining a complex mathematical structure is to sam\-ple i.i.d. random
points and look at some 
form of induced substructure relating the random points 
\end{quote}
which  assumes we are \emph{given} the complex structure.  
But now the second and more substantial part of the theme 
is that we can often  use exchangeability in the 
\emph{construction} of complex random structures as the $n \to \infty$ limits of 
random finite $n$-element  structures $\GG(n)$.  
\begin{quote}
Within the $n$-element structure $\GG(n)$ pick $k$
 random elements, look at the induced substructure on these $k$ elements---call this $\SS(n,k)$.
 Take a limit  (in distribution) as $n \to \infty$ for fixed $k$, any
necessary rescaling having been already done in the definition of
$\SS(n,k)$---call this limit $\SS_k$.
Within the limit random structures $(\SS_k, 2 \leq k < \infty)$, the $k$ elements
are exchangeable, and the distributions are consistent as $k$ increases and
therefore can be used to \emph{define} an infinite structure  $\SS_\infty$. 
\end{quote}
Where one can implement this program, the random structure $\SS_\infty$  will for
many purposes serve as a $n \to \infty$ limit of the original $n$-element
structures.  
Note that $\SS_\infty$ makes sense as a rather abstract object, via the
Kolmogorov\index{Kolmogorov, A. N.!Kolmogorov extension theorem}
extension theorem, but in concrete cases one tries to identify it with some more
concrete construction.

\subsection{First examples} 
To invent a name for the program above, let's call it 
\emph{exchangeable representations of complex random structures}\index{exchangeability!exchangeable representation}.
Let me first mention three examples.

\begin{enumerate}[{\bf 1.}]
\item Our discussion (section \ref{sec-partitions}) of exchangeable random partitions 
 fits the program  but is atypically simple, in that the 
limit $\SS_\infty$  is visibly the same kind of object (an \index{random partition}exchangeable random partition) as 
is the finite object $\GG(n)$.  
But when we moved on to  coalescence\index{coalescence} and
fragmentation\index{fragmentation} \emph{processes}  in section \ref{sec-bertoin},
our `exchangeability'  viewpoint prompts consideration of the limit process as the partition-valued 
process $(\Pi(t))$, which is rather different from the
finite-state \index{Markov, A. A.!Markov process|)}Markov processes arising in coalescence for finite $n$.

\item A conceptually similar example arises in the technically more
sophisticated setting of measure-valued
diffusions\index{diffusion!measure-valued diffusion} $(\mu_t)$.  
In such processes the  states are probability measures $\mu$ on some type-space, representing a `continuous'  infinite population.  
But one can alternately represent $\mu$  via an infinite i.i.d. sequence of samples from $\mu$, and thereby represent  the state more directly  as a 
discrete countable infinite population $(Z(i), i \geq 1)$ and the process as a \Index{particle process} $(Z_t(i), i \geq 1)$.
This viewpoint 
 was emphasized in the \emph{look-down construction}\index{look-down} of
Kurtz--Donnelly\index{Kurtz, T. G.}\index{Donnelly, P. [Donnelly, P. J.]}
\cite{MR1404525,MR1681126}.

\item As a complement to the characterization 
(\ref{vershik-cor}) of metric spaces with a probability measure (p.m.)\index{isometry}, we can define a notion of {\em
convergence} of such objects, say of finite metric spaces with p.m.s 
$(S_n,d_n,\mu_n)$ to a continuous limit $(S_\infty,d_\infty,\mu_\infty)$. The definition is simply that 
we have weak convergence 
\[ \bX^n \cd \bX^\infty \]
of the induced random arrays defined  as at (\ref{Xdxi}): 
\[ X^n_{\{i,j\}} = d_n(\xi^n_i, \xi^n_j);\  \{i,j\} \in \UP \]
for i.i.d. $(\xi^n_i, i \geq 1)$ with distribution $\mu_n$.
This definition provides an intriguing  complement to the more familiar notion
of \index{Gromov, M. L.!Gromov--Hausdorff distance}Gromov--Hausdorff distance between compact
metric spaces.
\end{enumerate}

Let us move on to the fundamental setting where the program gives a substantial payoff, the
formalization of  \emph{continuum random trees}\index{random tree|(} as rescaled limits of finite random
trees.

\subsection{Continuum random trees}
\label{sec-crt}
The material here is from the old survey \cite{me55}.

Probabilists are familiar with the notion that rescaled
\index{random walk (RW)}random walk converges in
distribution to \index{Brown, R.!Brownian motion}Brownian motion.
Now in the most basic case---simple symmetric RW of length $n$---we are studying 
the uniform distribution 
on a combinatorial set, the set of all $2^n$ simple walks of length $n$. 
So what happens if we study instead the uniform distribution on some other
combinatorial set? 
Let us consider the set of $n$-vertex trees. 
More precisely, consider either the set of rooted labeled trees
(Cayley's\index{Cayley, A.!Cayley's formula} formula
says there are 
$n^{n-1}$ such trees), or the set of rooted ordered trees (counted by the
Catalan\index{Catalan, E. C.!Catalan number}
number 
$\frac{1}{n} \binom{2n-2}{n-1}$), and write 
$\TT_n$ for the uniform random tree.

\begin{figure}[!ht]
\setlength{\unitlength}{0.33in}
\begin{picture}(10,6)(-1.2,-3)
\put(3,0.5){\line(-1,0){1.4}}
\put(3,0.5){\line(0,1){1.2}}
\put(2,-1){\line(-1,0){2}}
\put(2,-1){\line(0,-1){1}}
\put(4,0){\line(-2,1){1}}
\put(4,0){\line(-2,-1){2}}
\put(4,0){\line(1,0){5.5}}
\put(6,0){\line(0,-1){1.5}}
\put(7.3,0){\line(0,1){2.2}}
\put(0.93,0.32){\fbox{3}}
\put(-0.59,-1.18){\fbox{6}}
\put(1.72,-2.55){\fbox{2}}
\put(2.79,1.92){\fbox{1}}
\put(5.72,-2.05){\fbox{7}}
\put(7.07,2.43){\fbox{5}}
\put(9.59,-0.17){\fbox{4}}
\end{picture}
\caption{A leaf-labeled tree with edge-lengths.}\label{tree}
\end{figure}
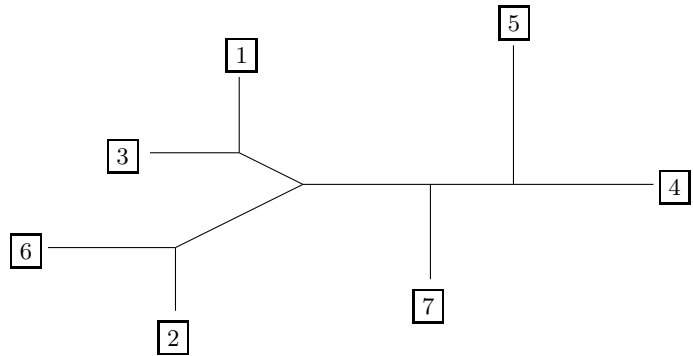

Trees fit nicely into the `substructure' framework. 
Vertices $v(1)$, \ldots, $v(k)$ of a tree define a spanning (sub)tree.  Take each
maximal path
$(w_0,w_1,\ldots,w_\ell)$ in the spanning tree whose intermediate vertices have 
degree $2$, and contract to a single edge of length $\ell$.
Applying this to $k$ independent uniform random vertices from a $n$-vertex model
$\TT_n$, 
then rescaling edge-lengths by the factor $n^{-1/2}$, gives a tree we'll call 
$\SS(n,k)$.  
We visualize such trees as in Figure \ref{tree}, vertex $v(i)$ having been relabeled as $i$.

\noindent
In the two models for $\TT_n$ mentioned above, one can do explicit calculations of
the distribution of 
$\SS(n,k)$, and use these to show that in distribution there is an $n \to \infty$
                             limit $\SS_k$ which
(up to a model-dependent scaling constant we'll ignore in this informal exposition) 
is the following distribution.

\begin{enumerate}[(i)]
\item The state space is the space of trees with $k$ leaves labeled 1,
2, \ldots, $k$ 
and with unlabeled degree-$3$ internal vertices, 
and where the $2k-3$ edge-lengths are positive real numbers. 

\item For each possible topological shape, the chance that the tree 
has that particular shape and that the vector of edge-lengths
$(L_1,\ldots,\allowbreak L_{2k-3})$ is in $([l_i, l_i + dl_i], 1 \leq i \leq 2k-3)$ 
equals $s \exp(-s^2/2) dl_1 \ldots\allowbreak dl_{2k-3}$, where $s = \sum_i l_i$. 
\end{enumerate}

One can check from the explicit formula what must be true from the general program,
that for fixed $k$ 
the distribution is exchangeable (in labels 1, \ldots, $k$),
and the distributions are consistent as $k$ 
increases (that is, the subtree of $\SS_{k+1}$ spanned by leaves 1, \ldots, $k$ is
distributed as $\SS_k$).

So some object $\SS_\infty$ exists, abstractly---but what is it, more concretely?
A \emph{real tree} \index{tree!real tree}is a metric space with the `tree' property 
that between any two points $v$ and $w$ there is a unique path. 
This implicitly specifies a length measure $\lambda$ such that the metric distance
$d(v,w)$ equals the length measure of the set of points on the path from $v$ to $w$. 
When a real tree is equipped with a mass measure $\mu$ 
of total mass $1$, representing a method for picking a vertex at random, I call it a 
\emph{continuum tree}. 
We will consider random continuum trees---which I call continuum random trees or 
CRTs because it sounds better!---and the \index{portmanteau theorem}Portmanteau Theorem below envisages realizations of $\SS_\infty$ 
as being equipped with a mass measure.  

Returning to the $n$-vertex random tree models $\TT_n$, by assigning\break `mass' $1/n$ to
each vertex we obtain the analogous `mass measure' on the vertices, used for 
randomly sampling vertices.  
The next result combines existence, construction and convergence theorems.  
The careful reader will notice that some details in the statements have been omitted.

\paragraph{The Portmanteau Theorem \cite{me55,me56}} 

\begin{enumerate}[{\bf 1.}]
\item {\bf Law of spanning subtrees.} \index{spanning subtree}There exists a particular \emph{Brownian CRT}
which agrees with 
$\SS_\infty$ in the following sense.   
Take a realization of the Brownian CRT, then pick $k$ i.i.d. vertices from the mass
measure, 
and consider the spanning subtree on these $k$ vertices. The unconditional law of
this subtree is the law 
in (ii) above. 

\item {\bf Construction from Brownian excursion.}\index{Brown, R.!Brownian excursion}
Consider an excursion-type function $f:[0,1] \to [0,\infty)$ with $f(0) = f(1) = 0$
and 
$f(x)>0$ elsewhere. Use $f$ to define a continuum tree as follows. 
Define a pseudo-metric on $[0,1]$ by: $d(x,y) = f(x) + f(y) - 2 \min(f(u): x \leq u
\leq y)$, $x\leq y$. 
The continuum tree is the associated metric space, and the mass measure is the 
image of Lebesgue measure on $[0,1]$. Using this construction with standard Brownian
excursion (scaled by a factor $2$) 
gives the Brownian CRT. 

\item {\bf Line-breaking construction.}\index{stick breaking}
Cut the line $[0,\infty)$ into finite segments at the points of a non-homogeneous
\index{Poisson, S. D.!Poisson process|(}Poisson process of 
intensity $\lambda(x) = x$. Build a tree by inductively attaching a segment
$[x_i,x_{i+1}]$ to a 
uniform random point of the tree built from the earlier segments. The tree built
from the first 
$k-1$ segments has the law (ii) above. The metric space closure of the tree built
from the whole half-line 
is the Brownian CRT, where the mass measure is the a.s. weak limit of the empirical
law of the first 
$k$ cut-points. 

\item {\bf Weak limit of conditioned critical Galton--Watson branching processes and of uniform random trees.} 
Take a critical \index{Galton, F.!Galton--Watson process}Galton--Watson branching process where the offspring law has finite
non-zero variance, 
and condition on total population until \Index{extinction} being $n$. This gives a random
tree. 
Rescale edge-lengths to have length $n^{-1/2}$. Put mass $1/n$ on each vertex. 
In a certain sense that can be formalized, the $n \to \infty$ weak limit of these
random trees 
is the Brownian CRT (up to a scaling factor). 
This result includes as special cases the two combinatorial models $\TT_n$ described
above.
\end{enumerate}

\subsection{Complements to the continuum random tree}
\label{sec-compl}
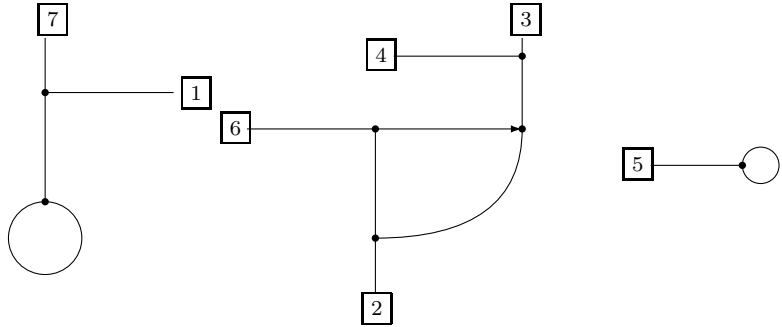
\begin{figure}[!ht]\footnotesize 
\setlength{\unitlength}{0.19in}
\begin{picture}(20,11.6)(-6,-6)
\put(-0.22,-0.15){\fbox{6}}
\put(7.7,2.82){\fbox{3}}
\put(3.75,1.85){\fbox{4}}
\put(3.63,-5.12){\fbox{2}}
\put(-5.2,2.82){\fbox{7}}
\put(-1.29,0.8){\fbox{1}}
\put(10.75,-1.15){\fbox{5}}
\put(0.5,0){\line(1,0){3.5}}
\put(4,0){\vector(1,0){4.0}}
\put(4,0){\line(0,-1){4.5}}
\put(8,0){\line(0,1){2.5}}
\put(8,2){\line(-1,0){3.5}}
\put(-5,-2){\line(0,1){4.5}}
\put(-5,1){\line(1,0){3.5}}
\put(-5,-3){\circle{2}}
\put(14,-1){\line(-1,0){2.5}}
\put(14.5,-1){\circle{1}}
\put(8,0){\circle*{0.2}}
\put(4,0){\circle*{0.2}}
\put(8,2){\circle*{0.2}}
\put(-5,-2){\circle*{0.2}}
\put(-5,1){\circle*{0.2}}
\put(14,-1){\circle*{0.2}}
\put(4,-3){\circle*{0.2}}
\qbezier(8,0)(8,-3)(4,-3)
\end{picture}
\caption{{\small The substructure $\SS_k$ of a random mapping.}}\label{randmap}
\end{figure}

More recent surveys by Le Gall\index{Le Gall, J.-F.} \cite{MR2203728} and by
Evans\index{Evans, S. N.} \cite{MR2351587} 
show different directions of development of the preceding material over the last 15 years.  
For instance
\begin{enumerate}[(i)]
\item  the \emph{Brownian snake}\index{Brown, R.!Brownian snake} \cite{MR2203728}, 
which combines the genealogical structure of random real trees with independent
spatial motions.
\item Diffusions on real trees: \cite{MR2351587} Chapter 7.
\item \emph{Continuum-tree valued
diffusions.}\index{diffusion!continuum-tree-valued diffusion}
There are several natural ways to define \index{Markov, A. A.!Markov chain}Markov chains on the space of $n$-vertex
trees 
such that the stationary distribution is uniform.  
Since the $n \to \infty$ rescaled limit of the stationary distribution is the Brownian CRT, 
it is natural to conjecture that the entire rescaled process can be made to converge
to some 
continuum-tree valued diffusion whose stationary distribution is the Brownian CRT.
But this forces us to engage a question that was deliberately avoided in the
previous section:
 what exactly is the space of all continuum trees, and when should we consider two
such trees to be the same?
This issue is discussed carefully in \cite{MR2351587}, based on the notion of the 
\index{Gromov, M. L.!Gromov--Hausdorff space}Gromov--Hausdorff space of all compact spaces.  Two specific continuum-tree valued
diffusions are then studied in 
Chapters 5 and 9 of \cite{MR2351587}.
\item Perhaps closer to our `exchangeability' focus, 
a surprising aspect of CRTs is their application to
stochastic \Index{coalescence}. 
For $0< \lambda < \infty$ split the Brownian CRT into components at the points of a 
\index{Poisson, S. D.!Poisson process|)}Poisson process of rate $\lambda$ along the skeleton of the tree. 
This gives a vector $Y(\lambda) = (Y_1(\lambda),Y_2(\lambda),\ldots)$ of masses 
of the components, which as $\lambda$ increases specifies a \index{fragmentation}fragmentation process. 
Reversing the direction of time by setting 
$\lambda = e^{-t}$ provides a construction of 
the (standard) additive coalescent \cite{me82}, 
that is the stochastic coalescent (section \ref{sec-bertoin}) with kernel $\kappa(x,y) = x+y$ 
`\Index{started from dust}'.
This result is  non-intuitive, and notable as one of  a handful of precise instances of the conceptual duality between stochastic  coalescence and fragmentation.
Also surprisingly, there are different ways that the additive coalescent can be `started from dust', and these can also be constructed via fragmentation of certain 
inhomogeneous CRTs \cite{me87}.   
This new family of CRTs satisfies analogs of the \index{portmanteau theorem}Portmanteau Theorem, 
and in particular there is an explicit analog of the formula (ii) in section \ref{sec-crt} 
for the distribution of the subtree $\SS_k$ spanned by $k$ random vertices \cite{me88}. 
This older work is complemented by much current work, the flavor of which can be seen in \cite{HPW09}.
\item A function $F: \{1,\ldots,n\} \to  \{1,\ldots,n\}$ defines a directed graph with edges $(i,F(i))$, and the topic 
\emph{random mappings}\index{random mapping} studies the graph derived  from a random function $F$.  
One can repeat the section \ref{sec-program} general program in this context. 
Any $k$ sampled vertices define an induced substructure, the subgraph of edges 
$i \to F(i) \to F(F(i)) \to \ldots$ reachable from some one of the sampled vertices.
Analogously to Figure \ref{tree}, contract paths between sampled vertices/junctions to single edges, to obtain (in the $n \to \infty$ limit) 
a \index{graph|(}graph $\SS_k$ with edge-lengths,
illustrated  in Figure \ref{randmap}. 
The theory of $n \to \infty$  limits of random mappings turns out to be closely related to that of random trees; the approach based on studying the 
consistent family $(\SS_k, k \ge 1)$ was developed in
Aldous--Pitman\index{Pitman, J. [Pitman, J. W.]} \cite{me98}.
\end{enumerate}\index{random tree|)}

\subsection{Second-order distance structure in random networks}

In the context (section \ref{sec-crt}) of continuum random trees the substructure was
distances between sampled points. 
At first sight one might hope that in many models of size-$n$ random
networks\index{random network|(} one
could repeat 
that analysis and find an interesting limit structure. 
But the particular feature of the models in section \ref{sec-crt} is `first order randomness' 
of the distance $D_n$ between two random vertices;
$D_n/ED_n$ has a non-constant limit distribution, leading to the randomness in the
limit structure.
Other models tend to fall into one of two categories.
For geometric networks (vertices having positions in $\Reals^2$; route-lengths as
Euclidean length) the route-length tends to grow as
constant $c$ times Euclidean distance, so any limit structure reflects only the
randomness of sampled vertex positions 
and the constant $c$, not any more interesting properties of the network. 
In non-geometric (e.g. Erd{\H{o}}s--R{\'e}nyi\index{Erd\H os,
P.}\index{Renyi, A.@R\'enyi, A.} random graph\index{graph|)}) models, $D_n$ tends to
be first-order constant.
So counter-intuitively, we don't know any other first-order random limit structures
outside the `somewhat tree-like' context.

Understanding \emph{second}-order behavior in spatial models is very chall\-enging---for instance, the second order behavior of first passage \Index{percolation} times remains a
longstanding open problem.
But one can get second order results in simple `random graph' type models, and here
is the basic example  
(mentioned in Aldous and Bhamidi\index{Bhamidi, S.} \cite{me119} as provable by the methods of that paper).
The probability model used for a random $n$-vertex network $\GG_n$ starts with the
complete graph
and assigns independent Exponential(rate $1/n$)
random lengths $L_{ij} = L_{ji} = L_e$ to the
$\binom n2$ edges $e = (i,j)$.
In this model 
$E D_n = \log n + O(1) $ 
and $\mathrm{var} (D_n) = O(1)$, and there is second-order behavior---a non-constant limit
distribution for 
$D_n - \log n$.

Now fix $k \geq 3$ and write $D_n(i,j) \ed D_n$ for the distance between vertices $i$ and $j$.
We expect a joint limit
\begin{equation}
(D_n(1,2) - \log n,\ldots,D_n(1,k) - \log n) \cd (D(1,2), \ldots, D(1,k))
\end{equation}
and it turns out the limit distribution is
\[ (D(1,2), \ldots, D(1,k)) \ed (\xi_1 + \eta_{12}, \ldots, \xi_1 + \eta_{1k}) \]
where $\xi_1$ has the double exponential distribution
\[ \Pr (\xi \leq x) = \exp(-e^{-x}), \ - \infty < x < \infty, \]
the $\eta_{1j}$ have \Index{logistic distribution}
\[ \Pr (\eta \leq x) = \sfrac{e^x}{1+e^x},  \ - \infty < x < \infty \]
and (here and below) the r.v.s in the limits are independent.
Now we can go one step further: we expect a joint limit for the array
\[ (D_n(i,j) - \log n, 1 \leq i < j \leq k) \cd (D(i,j), 1 \leq i < j \leq k) \]
and it turns out that the joint distribution of the limit is
\[ (D(i,j), 1 \leq i < j \leq k) \ed (\xi_i + \xi_j - \xi_{ij},  1 \leq i < j \leq
k)  \]
where the limit r.v.s all have the double exponential
distribution\index{Gumbel, E. J.!Gumbel distribution}.  
Of course the limit here must fit into the format of the partially exchangeable representation theorem (Theorem \ref{T1}),
and it is pleasant to see an explicit function $f$\index{random network|)}.

\section{Limits of finite deterministic structures}
Though we typically envisage limiting random structures arising as limits of finite \emph{random} structures, 
it also makes sense to consider limits of  finite \emph{deterministic} structures.  
Let me start with a trivial example. 
Suppose that for each $n$ we have a sequence $b_{n,1}$, \ldots, $b_{n,n}$ of $n$ bits (binary digits), and write $p_n$ for the proportion of $1$s.  
For each $k$  and $n$, sample $k$ random bits from $b_{n,1}$, \ldots,
$b_{n,n}$ and call the samples $X_{n,1}$, \ldots, $X_{n,k}$.  
Then, rather obviously, the property\undex{Bernoulli, J.}
\begin{align*}
 (X_{n,1},\ldots,X_{n,k}) &\cd \mbox{($k$ independent Bernoulli($p$))}\\
&{}\qquad\qquad \mbox{as $n \to \infty$;   for each $k$}
\end{align*}
is equivalent to the property $p_n \to p$.

Now in one sense this illustrates a big limitation to the  whole program\break---sampling  a substructure might lose most of the interesting information in the original structure!
But a parallel procedure in the deterministic graph setting (next section) does get more interesting results, and more sophisticated uses are mentioned in section \ref{sec-austin}.

\subsection{Limits of dense graphs}
\label{sec-dense}
Suppose that for each $n$ we have a \index{graph|(}graph $G_n$ on $n$ vertices.  
Write $p_n$ for the proportion of edges, relative to the total number $\binom n2$ of possible edges.  
We envisage the case $p_n \to p \in (0,1)$.  

For each $n$ let $(U_{n,i}, i \geq 1)$ be i.i.d. uniform on 1, \ldots, $n$.
Consider the infinite $\{0,1\}$-valued matrix $\bX^n$:
\[ X^n_{i,j} = 1((U_{n,i},U_{n,j}) \mbox{ is an edge of } G_n) . \]
When $n \gg k^2$ the $k$ sampled vertices $(U_{n,1}, \ldots, U_{n,k})$ of $G_n$ will
be distinct 
and the $k \times k$ restriction of $\bX^n$ is the incidence matrix of the 
induced subgraph $\SS(n,k)$ on these $k$ vertices.  
Suppose there is a limit \Index{random matrix} $\bX$:
\begin{equation}
\bX^n \cd \bX \mbox{ as } n \to \infty 
\label{bXX}
\end{equation}
in the usual \Index{product topology}, that is 
\[ (X^n_{i,j},\  1 \leq i,j \leq k) \cd (X_{i,j},\  1 \leq i,j \leq k) \mbox{ for
each } k . \]
(Note that by compactness there is always a \emph{subsequence} in which such
convergence holds.)
Now each $\bX^n$ has the partially exchangeable property (\ref{def-PE}), and the
limit $\bX$ 
 inherits this property, so we can apply  the representation
 theorem (Theorem \ref{T1}) to describe
the possible limits.
In the $\{0,1\}$-valued case we can simplify the representation.  First consider a 
function of form  (\ref{PE-rep}) but not depending on the first coordinate---that is, a function 
$f(u_i,u_j, u_{\{i,j\}}  )$.  Write 
\[ q(u_i,u_j) = \Pr (f(u_i,u_j, u_{\{i,j\}}  ) = 1).\]
The distribution of a $\{0,1\}$-valued \index{exchangeability!partial exchangeability|(}partially exchangeable array of the special
form 
$f(U_i,U_j, U_{\{i,j\}}  )$ is determined by the symmetric function $q(\cdot,\cdot)$, 
and so for the general form (\ref{PE-rep}) the distribution is specified by a probability
distribution over such 
symmetric functions.  

This all fits our section \ref{sec-program} general program.
From an arbitrary sequence of finite deterministic graphs we can (via passing to a
subsequence if necessary) 
extract a `limit infinite random graph' $\SS_\infty$ on vertices 1, 2, \ldots,
defined by its incidence matrix $\bX$ in the limit (\ref{bXX}), 
and we can characterize the possible limits.  
But what is a more concrete interpretation of the relation between $\SS_\infty$ and
the finite graphs $(G_n)$?
To a probabilist the verbal expression of (\ref{bXX})
\begin{quote}
the restriction $\SS_k$ of $\SS_\infty$ to vertices 1, \ldots, $k$ is distributed
as the $n \to \infty$ limit of the 
induced subgraph of $G_n$ on $k$ random vertices
\end{quote}
is clear enough, but here is a translation into more graph-theoretic language,
following \cite{MR2463439}.
For finite graphs $F$, $G$ write $t(F,G)$ for the proportion of all mappings from
vertices of $F$ to 
vertices of $G$ that are graph homomorphisms, i.e. map adjacent vertices to adjacent
vertices.  
Suppose $F$ has $k$ vertices, and we label them arbitrarily as 1, \ldots, $k$. 
Take the subgraph $G[k]$ of $G$ on $k$ randomly sampled vertices, labeled
1, \ldots, $k$, and note that whether we sample with or without
replacement\index{sampling} makes no difference to $n \to \infty$ limits.
Then $t(F,G)$ is the probability that $F$ is a subgraph of $G[k]$.
Now write $t_=(F,G)$ for the probability that $F = G[k]$.  
For fixed $k$, a standard \Index{inclusion-exclusion} argument shows that, for a sequence
$(G_n)$,
the existence of either family of limits
\begin{eqnarray}
 \lim_n t(F,G_n) \mbox{ exists, for each graph $F$ on vertices } \{1,\ldots,k\}, &&
\label{tFG} \\
 \lim_n t_=(F,G_n) \mbox{ exists, for each graph $F$ on vertices } \{1,\ldots,k\}, &&
\label{tFG=}
\end{eqnarray}
implies existence of the other family of limits.

In our program, the notion of $\SS_\infty$ being the limit of $G_n$ was defined by
(\ref{bXX}),
which is equivalent to requiring existence of limits (\ref{tFG=}) for each $k$, in
which case the limits are 
just $\Ex t_=(F,\SS_k)$.  And as indicated above, the  
partially exchangeable representation theorem (Theorem \ref{T1}) 
characterizes the possible limit structures $\SS_\infty$.
A recent line of work in graph theory, initiated by Lov{\'a}sz\index{Lov\'asz,
L.} and Szegedy\index{Szegedy, B.}
\cite{MR2274085},
started by defining convergence in the equivalent way via (\ref{tFG}) and obtained
the same characterization.
This is the second recent rediscovery of special cases of partially exchangeable
representation theory\index{exchangeability!partial exchangeability|)}.
Diaconis\index{Diaconis, P.} and Janson\index{Janson, S.} \cite{MR2463439} give a very clear and detailed account of the 
relation between the two settings,
and Diaconis--Holmes--Janson\index{Holmes, S.} \cite{persi-2009} work through to an explicit 
description of the possible limits for a particular subclass of graphs 
called \emph{threshold graphs}\index{graph!threshold graph}.
Of course the line of work started in \cite{MR2274085} has been developed further  to produce new and interesting results in graph theory---see e.g. \cite{MR2455626}\index{graph|)}.

\subsection{Further uses in finitary combinatorics}\index{combinatorics}
\label{sec-austin}
The remarkable recent survey by Austin\index{Austin, T.} \cite{MR2426176} gives a more sophisticated treatment of the 
theory of representations of jointly exchangeable arrays, with the goal (\cite{MR2426176}
section 4) 
of clarifying connections between that theory and topics involving limits in
finitary combinatorics, such as 
those in our previous section.  
I don't understand this material well enough to do more than copy a few phrases, as
follows. 
Section 4.1 of \cite{MR2426176} gives a general discussion of `extraction of limit objects',
somewhat parallel to our
section \ref {sec-program}, but with more detailed discussion of different possible precise
mathematical structures.  
The paper continues, describing connections with the `hypergraph regularity lemmas' 
 featuring in combinatorial proofs of Szemer\'{e}di's
Theorem\index{Szemer\'edi, E.!Szemer\'edi's theorem}, and with the structure
theory within \Index{ergodic theory}
that Furstenberg\index{Furstenberg, H.} developed for his proof of Szemer\'{e}di's Theorem.
A subsequent technical paper Austin--Tao\index{Tao, T.} \cite{austin-2008} applies such methods to
the topic 
of hereditary properties of  graphs or hypergraphs being testable with one-sided
error; 
informally, this means that if a graph or \Index{hypergraph} satisfies that property 
`locally' with sufficiently high probability, then it can be 
modified into a graph or hypergraph which satisfies that property `globally'.

\section{Miscellaneous comments}
\label{sec-misc}
\begin{enumerate}[{\bf 1.}]
\item To get an idea of the breadth of the topic,
\emph{Mathematical Reviews}\index{Mathematical Reviews@{\it Mathematical Reviews}} created an
`exchangeability' classification 60G09 in 1984, which has attracted around 300
items; 
\emph{Google Scholar}\index{Google@{\it Google}!Google Scholar@{\it Google Scholar}} finds around 350 citations of the survey \cite{me22}; and the
overlap is only around 50\%.  
The topics in this paper, centered around structure theory---theory and applications of extensions of de
Finetti's\index{Finetti, B. de!de Finetti's theorem} theorem---are in fact only a rather small part of this whole.  
In particular the 
`exchangeable pairs' idea central to Stein's\index{Stein, C. M.!Stein's method} method  \cite{MR2201882} is really a completely distinct field.

\item Our central theme involved exchangeability, but one can perhaps view it as part of a
broader theme:
\begin{quote} 
a mathematical object equipped with a probability measure  is sometimes a richer and
more natural structure than the object by itself.
\end{quote}
For instance, elementary discussion of fractals\index{fractal} like the
Sierpi\'nski\index{Sierpinski, W.@Sierpi\'nski, W.!Sierpi\'nski gasket} gasket view the
object as a set in $\Reals^2$, but it comes equipped with its natural `uniform
probability distribution' which enables richer questions---the measure of small balls  around a
typical point, for example.  
Weierstrass's\index{Weierstrass, K. T. W.} construction of a continuous \Index{nowhere differentiable} function seems at
first sight artificial---where would such things arise naturally?---but then
the fact that the \index{Brown, R.!Brownian motion}Brownian motion process puts a probability measure on such functions 
indicates one place where they do arise naturally.    Analogously the notion of \emph{real
tree}\index{tree!real tree} (section \ref{sec-crt}) may seem at first sight artificial---how might such objects
arise naturally?---but then realizing they arise as limits of random finite trees 
indicates one place where they do arise naturally.  
 Of course the underlying structure `a space with a metric and a measure' arises 
in many contexts, for example (under the name 
\emph{metric measure space})\index{metric measure space} in the context of \Index{differential geometry} questions 
\cite{MR2408268}.
\end{enumerate}\index{exchangeability|)}

\paragraph{Acknowledgments}
I thank Persi Diaconis\index{Diaconis, P.}, Jim Pitman\index{Pitman, J. [Pitman, J. W.]} and an anon\-ymous referee for very helpful comments.

\bibliographystyle{cambridgeauthordateCMG}
 \bibliography{aldous_kingman}

\begin{thebibliography}{58}
\expandafter\ifx\csname natexlab\endcsname\relax\def\natexlab#1{#1}\fi
\expandafter\ifx\csname selectlanguage\endcsname\relax
  \def\selectlanguage#1{\relax}\fi

\bibitem[\protect\citename{Aldous, }1977]{me2}
Aldous, D.~J. 1977.
\newblock Limit theorems for subsequences of arbitrarily-dependent sequences of
  random variables.
\newblock {\em Z. Wahrscheinlichkeitstheorie verw. Gebiete}, {\bf 40}(1),
  59--82.

\bibitem[\protect\citename{Aldous, }1981]{me11}
Aldous, D.~J. 1981.
\newblock Representations for partially exchangeable arrays of random
  variables.
\newblock {\em J. Multivariate Anal.}, {\bf 11}, 581--598.

\bibitem[\protect\citename{Aldous, }1985]{me22}
Aldous, D.~J. 1985.
\newblock Exchangeability and related topics.
\newblock {Pages  1--198 of:} Hennequin, P.~L. (ed), {\em {\'E}cole
  d'{\'E}t{\'e} de Saint-Flour XIII---1983}.
\newblock Lecture Notes in Math., vol. 1117.
\newblock Berlin: Springer-Verlag.

\bibitem[\protect\citename{Aldous, }1991]{me55}
Aldous, D.~J. 1991.
\newblock The continuum random tree {II}: an overview.
\newblock {Pages  23--70 of:} Barlow, M.~T., and Bingham, N.~H. (eds), {\em
  Stochastic Analysis}.
\newblock Cambridge: Cambridge Univ. Press.

\bibitem[\protect\citename{Aldous, }1993]{me56}
Aldous, D.~J. 1993.
\newblock The continuum random tree {III}.
\newblock {\em Ann. Probab.}, {\bf 21}, 248--289.

\bibitem[\protect\citename{Aldous, }1999]{me78}
Aldous, D.~J. 1999.
\newblock Deterministic and stochastic models for coalescence (aggregation and
  coagulation): a review of the mean-field theory for probabilists.
\newblock {\em Bernoulli}, {\bf 5}, 3--48.

\bibitem[\protect\citename{Aldous and Bhamidi, }2007]{me119}
Aldous, D.~J., and Bhamidi, S. 2007.
\newblock {\em Edge Flows in the Complete Random-Lengths Network}.
\newblock \url{http://arxiv.org/abs/0708.0555}.

\bibitem[\protect\citename{Aldous and Pitman, }1998]{me82}
Aldous, D.~J., and Pitman, J. 1998.
\newblock The standard additive coalescent.
\newblock {\em Ann. Probab.}, {\bf 26}, 1703--1726.

\bibitem[\protect\citename{Aldous and Pitman, }1999]{me88}
Aldous, D.~J., and Pitman, J. 1999.
\newblock A family of random trees with random edge lengths.
\newblock {\em Random Structures \& Algorithms}, {\bf 15}, 176--195.

\bibitem[\protect\citename{Aldous and Pitman, }2000]{me87}
Aldous, D.~J., and Pitman, J. 2000.
\newblock Inhomogeneous continuum random trees and the entrance boundary of the
  additive coalescent.
\newblock {\em Probab. Theory Related Fields}, {\bf 118}, 455--482.

\bibitem[\protect\citename{Aldous and Pitman, }2002]{me98}
Aldous, D.~J., and Pitman, J. 2002.
\newblock Invariance principles for non-uniform random mappings and trees.
\newblock {Pages  113--147 of:} Malyshev, V.~A., and Vershik, A.~M. (eds), {\em
  Asymptotic Combinatorics with Applications to Mathematical Physics}.
\newblock Dordrecht: Kluwer.
\newblock Available at
  \url{http://www.stat.berkeley.edu/tech-reports/594.ps.Z}.

\bibitem[\protect\citename{Austin, }2008]{MR2426176}
Austin, T. 2008.
\newblock On exchangeable random variables and the statistics of large graphs
  and hypergraphs.
\newblock {\em Probab. Surv.}, {\bf 5}, 80--145.

\bibitem[\protect\citename{Austin and Tao, }2008]{austin-2008}
Austin, T., and Tao, T. 2008.
\newblock {\em On the Testability and Repair of Hereditary Hypergraph
  Properties}.
\newblock \url{http://arxiv.org/abs/0801.2179}.

\bibitem[\protect\citename{Barbour and Chen, }2005]{MR2201882}
Barbour, A.~D., and Chen, L. H.~Y. (eds). 2005.
\newblock {\em Stein's Method and Applications}.
\newblock Lecture Notes Series, vol. 5.
\newblock Singapore: Singapore Univ. Press, for Institute for Mathematical
  Sciences, National University of Singapore.

\bibitem[\protect\citename{Berkes and P{\'e}ter, }1986]{MR859840}
Berkes, I., and P{\'e}ter, E. 1986.
\newblock Exchangeable random variables and the subsequence principle.
\newblock {\em Probab. Theory Related Fields}, {\bf 73}(3), 395--413.

\bibitem[\protect\citename{Bertoin, }2006]{MR2253162}
Bertoin, J. 2006.
\newblock {\em Random Fragmentation and Coagulation Processes}.
\newblock Cambridge Stud. Adv. Math., vol. 102.
\newblock Cambridge: Cambridge Univ. Press.

\bibitem[\protect\citename{Billingsley, }1968]{MR0233396}
Billingsley, P. 1968.
\newblock {\em Convergence of Probability Measures}.
\newblock New York: John Wiley \& Sons.

\bibitem[\protect\citename{Billingsley, }1995]{MR1324786}
Billingsley, P. 1995.
\newblock {\em Probability and Measure}. Third edn.
\newblock Wiley Ser. Probab. Math. Stat.
\newblock New York: John Wiley \& Sons.

\bibitem[\protect\citename{Bingham {et~al.}, }2003]{MR2026570}
Bingham, N.~H., Kiesel, R., and Schmidt, R. 2003.
\newblock A semi-parametric approach to risk management.
\newblock {\em Quant. Finance}, {\bf 3}(6), 426--441.

\bibitem[\protect\citename{Blei {et~al.}, }2003]{LDA2003}
Blei, D.~M., Ng, A.~Y., and Jordan, M.~I. 2003.
\newblock Latent {D}irichlet allocation.
\newblock {\em J. Machine Learning Res.}, {\bf 3}, 993--1022.

\bibitem[\protect\citename{Borgs {et~al.}, }2008]{MR2455626}
Borgs, C., Chayes, J.~T., Lov{\'a}sz, L., S{\'o}s, V.~T., and Vesztergombi, K.
  2008.
\newblock Convergent sequences of dense graphs. {I}. {S}ubgraph frequencies,
  metric properties and testing.
\newblock {\em Adv. Math.}, {\bf 219}(6), 1801--1851.

\bibitem[\protect\citename{Chow and Teicher, }1997]{MR1476912}
Chow, Y.~S., and Teicher, H. 1997.
\newblock {\em Probability Theory: Independence, Interchangeability,
  Martingales}. Third edn.
\newblock Springer Texts Statist.
\newblock New York: Springer-Verlag.

\bibitem[\protect\citename{Diaconis and Freedman, }1980a]{MR556418}
Diaconis, P., and Freedman, D.~A. 1980a.
\newblock de {F}inetti's theorem for {M}arkov chains.
\newblock {\em Ann. Probab.}, {\bf 8}(1), 115--130.

\bibitem[\protect\citename{Diaconis and Freedman, }1980b]{MR577313}
Diaconis, P., and Freedman, D.~A. 1980b.
\newblock Finite exchangeable sequences.
\newblock {\em Ann. Probab.}, {\bf 8}(4), 745--764.

\bibitem[\protect\citename{Diaconis and Freedman, }1984]{MR786142}
Diaconis, P., and Freedman, D.~A. 1984.
\newblock Partial exchangeability and sufficiency.
\newblock {Pages  205--236 of:} {\em Statistics: Applications and New
  Directions ({C}alcutta, 1981)}.
\newblock Calcutta: Indian Statist. Inst.

\bibitem[\protect\citename{Diaconis and Janson, }2008]{MR2463439}
Diaconis, P., and Janson, S. 2008.
\newblock Graph limits and exchangeable random graphs.
\newblock {\em Rend. Mat. Appl. (7)}, {\bf 28}(1), 33--61.

\bibitem[\protect\citename{Diaconis {et~al.}, }2009]{persi-2009}
Diaconis, P., Holmes, S., and Janson, S. 2009.
\newblock {\em Threshold Graph Limits and Random Threshold Graphs}.
\newblock \url{http://arxiv.org/abs/0908.2448}.

\bibitem[\protect\citename{Donnelly and Kurtz, }1996]{MR1404525}
Donnelly, P., and Kurtz, T.~G. 1996.
\newblock A countable representation of the {F}leming-{V}iot measure-valued
  diffusion.
\newblock {\em Ann. Probab.}, {\bf 24}(2), 698--742.

\bibitem[\protect\citename{Donnelly and Kurtz, }1999]{MR1681126}
Donnelly, P., and Kurtz, T.~G. 1999.
\newblock Particle representations for measure-valued population models.
\newblock {\em Ann. Probab.}, {\bf 27}(1), 166--205.

\bibitem[\protect\citename{Dovbysh and Sudakov, }1982]{MR666087}
Dovbysh, L.~N., and Sudakov, V.~N. 1982.
\newblock Gram-de {F}inetti matrices.
\newblock {\em Zap. Nauchn. Sem. Leningrad. Otdel. Mat. Inst. Steklov. (LOMI)},
  {\bf 119}, 77--86, 238, 244--245.
\newblock Problems of the Theory of Probability Distribution, {VII}.

\bibitem[\protect\citename{Durrett, }2005]{dur91v3}
Durrett, R. 2005.
\newblock {\em Probability: Theory and Examples}. Third edn.
\newblock Statistics/Probability Series.
\newblock Pacific Grove, CA: Brooks/Cole.

\bibitem[\protect\citename{Evans, }2008]{MR2351587}
Evans, S.~N. 2008.
\newblock {\em Probability and Real Trees}.
\newblock Lecture Notes in Math., vol. 1920.
\newblock Berlin: Springer-Verlag.

\bibitem[\protect\citename{Ewens and Watterson, }2009]{EW10}
Ewens, W.~J., and Watterson, G.~A. 2009.
\newblock Kingman and mathematical population genetics.
\newblock {In:} Bingham, N.~H., and Goldie, C.~M. (eds), {\em Probability and
  Mathematical Genetics: Papers in Honour of {S}ir {J}ohn {K}ingman}.
\newblock London Math. Soc. Lecture Note Ser.
\newblock Cambridge: Cambridge Univ. Press.

\bibitem[\protect\citename{Freedman, }1962]{MR0156369}
Freedman, D.~A. 1962.
\newblock Invariants under mixing which generalize de {F}inetti's theorem.
\newblock {\em Ann. Math. Statist}, {\bf 33}, 916--923.

\bibitem[\protect\citename{Fristedt and Gray, }1997]{MR1422917}
Fristedt, B., and Gray, L. 1997.
\newblock {\em A Modern Approach to Probability Theory}.
\newblock Probab. Appl. Ser.
\newblock Boston, MA: Birkh\"auser.

\bibitem[\protect\citename{Gelman {et~al.}, }2004]{MR2027492}
Gelman, A., Carlin, J.~B., Stern, H.~S., and Rubin, D.~B. 2004.
\newblock {\em Bayesian Data Analysis}. Second edn.
\newblock Texts Statist. Sci. Ser.
\newblock Boca Raton, FL: Chapman \& Hall/CRC.

\bibitem[\protect\citename{Gnedin {et~al.}, }2009]{GHP09}
Gnedin, A.~V., Haulk, C., and Pitman, J. 2009.
\newblock Characterizations of exchangeable partitions and random discrete
  distributions by deletion properties.
\newblock {In:} Bingham, N.~H., and Goldie, C.~M. (eds), {\em Probability and
  Mathematical Genetics: Papers in Honour of {S}ir {J}ohn {K}ingman}.
\newblock London Math. Soc. Lecture Note Ser.
\newblock Cambridge: Cambridge Univ. Press.

\bibitem[\protect\citename{Haas {et~al.}, }2009]{HPW09}
Haas, B., Pitman, J., and Winkel, M. 2009.
\newblock Spinal partitions and invariance under re-rooting of continuum random
  trees.
\newblock {\em Ann. Probab.}, {\bf 37}(4), 1381--1411.

\bibitem[\protect\citename{Hoover, }1979]{hoover-rel}
Hoover, D.~N. 1979.
\newblock {\em Relations on Probability Spaces and Arrays of Random Variables}.
\newblock Preprint, Institute of Advanced Studies, Princeton.

\bibitem[\protect\citename{Kallenberg, }1973]{MR0394842}
Kallenberg, O. 1973.
\newblock Canonical representations and convergence criteria for processes with
  interchangeable increments.
\newblock {\em Z. Wahrscheinlichkeitstheorie verw. Gebiete}, {\bf 27}, 23--36.

\bibitem[\protect\citename{Kallenberg, }1989]{MR1003713}
Kallenberg, O. 1989.
\newblock On the representation theorem for exchangeable arrays.
\newblock {\em J. Multivariate Anal.}, {\bf 30}(1), 137--154.

\bibitem[\protect\citename{Kallenberg, }1997]{kall97}
Kallenberg, O. 1997.
\newblock {\em Foundations of Modern Probability}.
\newblock Probab. Appl. (N.Y.).
\newblock New York: Springer-Verlag.

\bibitem[\protect\citename{Kallenberg, }2005]{MR2161313}
Kallenberg, O. 2005.
\newblock {\em Probabilistic Symmetries and Invariance Principles}.
\newblock Probab. Appl. (N.Y.).
\newblock New York: Springer-Verlag.

\bibitem[\protect\citename{Kingman, }1978]{MR0494344}
Kingman, J. F.~C. 1978.
\newblock Uses of exchangeability.
\newblock {\em Ann. Probab.}, {\bf 6}(2), 183--197.

\bibitem[\protect\citename{Le~Gall, }2005]{MR2203728}
Le~Gall, J.-F. 2005.
\newblock Random trees and applications.
\newblock {\em Probab. Surv.}, {\bf 2}, 245--311 (electronic).

\bibitem[\protect\citename{Lo{\`e}ve, }1978]{MR0651018}
Lo{\`e}ve, M. 1978.
\newblock {\em Probability Theory {II}}. Fourth edn.
\newblock Grad. Texts in Math., vol. 46.
\newblock New York: Springer-Verlag.

\bibitem[\protect\citename{Lott, }2007]{MR2408268}
Lott, J. 2007.
\newblock Optimal transport and {R}icci curvature for metric-measure spaces.
\newblock {Pages  229--257 of:} {\em Surveys in Differential Geometry}.
\newblock Surv. Differ. Geom., vol. 11.
\newblock Somerville, MA: Int. Press.

\bibitem[\protect\citename{Lov{\'a}sz and Szegedy, }2006]{MR2274085}
Lov{\'a}sz, L., and Szegedy, B. 2006.
\newblock Limits of dense graph sequences.
\newblock {\em J. Combin. Theory Ser. B}, {\bf 96}(6), 933--957.

\bibitem[\protect\citename{Nualart, }2006]{MR2200233}
Nualart, D. 2006.
\newblock {\em The {M}alliavin Calculus and Related Topics}. Second edn.
\newblock Probab. Appl. Ser.
\newblock Berlin: Springer-Verlag.

\bibitem[\protect\citename{Olshanski and Vershik, }1996]{MR1402920}
Olshanski, G., and Vershik, A. 1996.
\newblock Ergodic unitarily invariant measures on the space of infinite
  {H}ermitian matrices.
\newblock {Pages  137--175 of:} {\em Contemporary Mathematical Physics}.
\newblock Amer. Math. Soc. Transl. Ser. 2, vol. 175.
\newblock Providence, RI: Amer. Math. Soc.

\bibitem[\protect\citename{Panchenko, }2009]{panch09}
Panchenko, D. 2009.
\newblock {\em On the {D}ovbysh-{S}udakov Representation Result}.
\newblock \url{http://front.math.ucdavis.edu/0905.1524}.

\bibitem[\protect\citename{Pitman, }2006]{MR2245368}
Pitman, J. 2006.
\newblock {\em Combinatorial Stochastic Processes}.
\newblock Lecture Notes in Math., vol. 1875.
\newblock Berlin: Springer-Verlag.

\bibitem[\protect\citename{Pitman and Yor, }1997]{MR1434129}
Pitman, J., and Yor, M. 1997.
\newblock The two-parameter {P}oisson-{D}irichlet distribution derived from a
  stable subordinator.
\newblock {\em Ann. Probab.}, {\bf 25}(2), 855--900.

\bibitem[\protect\citename{Schoenberg, }1938]{MR1501980}
Schoenberg, I.~J. 1938.
\newblock Metric spaces and positive definite functions.
\newblock {\em Trans. Amer. Math. Soc.}, {\bf 44}(3), 522--536.

\bibitem[\protect\citename{{Smoluchowski}, }1916]{smoluch}
{Smoluchowski}, M.~von. 1916.
\newblock Drei {V}ortr{\"a}ge {\"u}ber {D}iffusion, {B}rownsche {B}ewegung und
  {K}oagulation von {K}olloidteilchen.
\newblock {\em Physik. Z.}, {\bf 17}, 557--585.

\bibitem[\protect\citename{Tak{\'a}cs, }1967]{MR0217858}
Tak{\'a}cs, L. 1967.
\newblock {\em Combinatorial Methods in the Theory of Stochastic Processes}.
\newblock New York: John Wiley \& Sons.

\bibitem[\protect\citename{Vershik, }2004]{MR2086637}
Vershik, A.~M. 2004.
\newblock Random metric spaces and universality.
\newblock {\em Uspekhi Mat. Nauk}, {\bf 59}(2(356)), 65--104.

\bibitem[\protect\citename{Wakeley, }2008]{wakeley}
Wakeley, J. 2008.
\newblock {\em Coalescent Theory: An Introduction}.
\newblock Greenwood Village, CO: Roberts \& Co.

\end{thebibliography}

\end{document}